\RequirePackage{fix-cm}
\documentclass[twocolumn]{svjour3}          
\smartqed  
\usepackage{graphicx,float}
\usepackage{amsmath,amsfonts,latexsym,amssymb}

\begin{document}

\title{Swimming by switching
}


\author{Fabio Bagagiolo        \and Rosario Maggistro \and Marta Zoppello
}


\institute{F. Bagagiolo\at
              Dipartimento di Matematica, Universit\'a degli studi di Trento \\
              via Sommarive,14 38123 Povo (TN)
              \email{fabio.bagagiolo@unitn.it}           
           \and
           M. Zoppello \at
            Dipartimento di Matematica, Universit\'a degli studi di Trento \\
             via Sommarive,14 38123 Povo (TN)
               \email{marta.zoppello@unitn.it}
               \and
             R. Maggistro\at
              Dipartimento di Matematica, Universit\'a degli studi di Trento \\
             via Sommarive,14 38123 Povo (TN)
             \email{rosario.maggistro@unitn.it}
}

\date{Received: date / Accepted: date}

\maketitle

\begin{abstract}
In this paper we investigate different strategies to overcome the scallop theorem. We will show how to obtain a net motion exploiting the fluid's type change during a periodic deformation. We are interested in two different models: in the first one that change is linked to the magnitude of the opening and closing velocity. Instead, in the second one it is related to the sign of the above velocity.
An interesting feature of the latter model is the introduction of a delay-switching rule through a thermostat. We remark that the latter is fundamental in order to get both forward and backward motion.

\end{abstract}

\section{Introduction}
\label{intro}
The study of locomotion strategies in fluids is attracting increasing interest in recent literature, especially for its connection with the realization of artificial devices that can self-propel in fluids. Theories of swimming generally utilize either low Reynolds number approximation, or the assumption of inviscid ideal fluid dynamics (high Reynolds number). These two different regimes are also distinct in terms of the mechanism of locomotion \cite{Childress81,Lighthill75}.\\ 
In this paper we focus on swimmers immersed in these two kind of fluids which produce a linear dynamics. In particular we study the system describing the motion of a scallop for which it is well known \cite{AlougesDeSimone08,GeneralizedScallop,Purcell77} that the \textit{scallop theorem/paradox} holds. This means that it is not capable to achieve any net motion performing cyclical shape changes, either in a viscous or in an inviscid fluid. Some authors tried to overcome this paradox changing the geometry of the swimmer, for example adding a degree of freedom, introducing the Purcell swimmer \cite{Purcell77}, or the three sphere swimmer \cite{Goldstein}. Others, instead, supposed the scallop immersed in a non Newtonian fluid, in which the viscosity is not constant, ending up with a non reversible dynamics \cite{Pseudoelastic,Nature}. Inspired by this last approach, our aim is to propose some strategies which maintain the swimmer geometry and exploit instead a change in the dynamics. The idea is based on switching dynamics depending on the angular velocity of opening and closing of the scallop's valves. More precisely we analyze two cases: in the first one we suppose that if the modulus of the angular velocity is high, the fluid regime can be approximated by the ideal one, instead if this modulus is low the fluid can be considered as completely viscous. These assumptions are realistic since the Reynolds number changes depending on the characteristic velocity of the swimmer. In the second case we assume that the fluid reacts in a different way between the opening and closing of the valves: it facilitates the opening, so that it can be considered an ideal fluid, and resists the closing, like a viscous fluid. These last approximations model a fluid in which the viscosity changes with the sign of the angular velocity. More precisely we use two constant viscosities: one high (resp. one very small) if the angular velocity is negative (resp. positive).  Moreover inspired by \cite{Nature}, where the scallop's opening and closing is actuated by an external magnetic field, in this last case we also introduce an hysteresis mechanism through a thermostat, see Fig \ref{Fig00} (see \cite{Visintin} for mathematical models for hysteresis), to model a delay in the change of fluid's regime. In both cases we assume to be able to prescribe the angular velocity, using it as a control parameter and we prove that the system is controllable, i.e. the scallop is able to move both forward and backward using cyclical deformations. Furthermore we prove also that it is always possible to move between two fixed points, starting and ending with two prescribed angles.\\ In the last part of the paper we show also some numerical examples to support our theoretical predictions.\\\\
The plan of the paper is the following. In Section \ref{sec:2} we present the swimmer model and derive its equation of motion both in the viscous and in the ideal approximation, proving the scallop theorem. Section \ref{sec:3} is devoted to the introduction of the switching strategies which lead to the controllability of the scallop system. Finally in Section \ref{sec:4} we present some numerical simulations showing different kind of controls that can be used.
\section{The Scallop swimmer}
\label{sec:2}
In this section we are interested in analyzing the motion of an articulated rigid body immersed in a fluid that changes its configuration. In order to determine completely its state we need the position of its center of mass and its orientation. Their temporal evolution is obtained solving the Newton's equations 
coupled with the Navier-Stokes equations relative to the surrounding fluid. We will face this problem considering the body as immersed in two kinds of different fluids: one viscous at low Reynolds number in which we neglect the effects of inertia, 
and another one ideal inviscid and irrotational, in which we neglect the viscous forces in the Navier-Stokes equations.  
First of all we recall that in both cases a swimmer that tries to moves like a scallop, opening and closing periodically its valves, does not move at the end of a cycle. This situation is well known as scallop theorem (or paradox) \cite{AlougesDeSimone08,Purcell77}. \\
In what follows we will consider a planar  body composed by two rigid valves of elliptical shape, joined in order that they can be opened and closed. Moreover this body is constrained to move along one of the cartesian axes (the $\vec{e}_x$-axis) and is symmetric with respect to it. Finally we will neglect the interaction between the two valves. The configuration of the system is easily described by the position $x$ of the juncture point along the  $\vec{e}_x$-axis and by the angle $\theta$ that each valve forms with the axis

\begin{figure}[H]
  \includegraphics[scale=0.35]{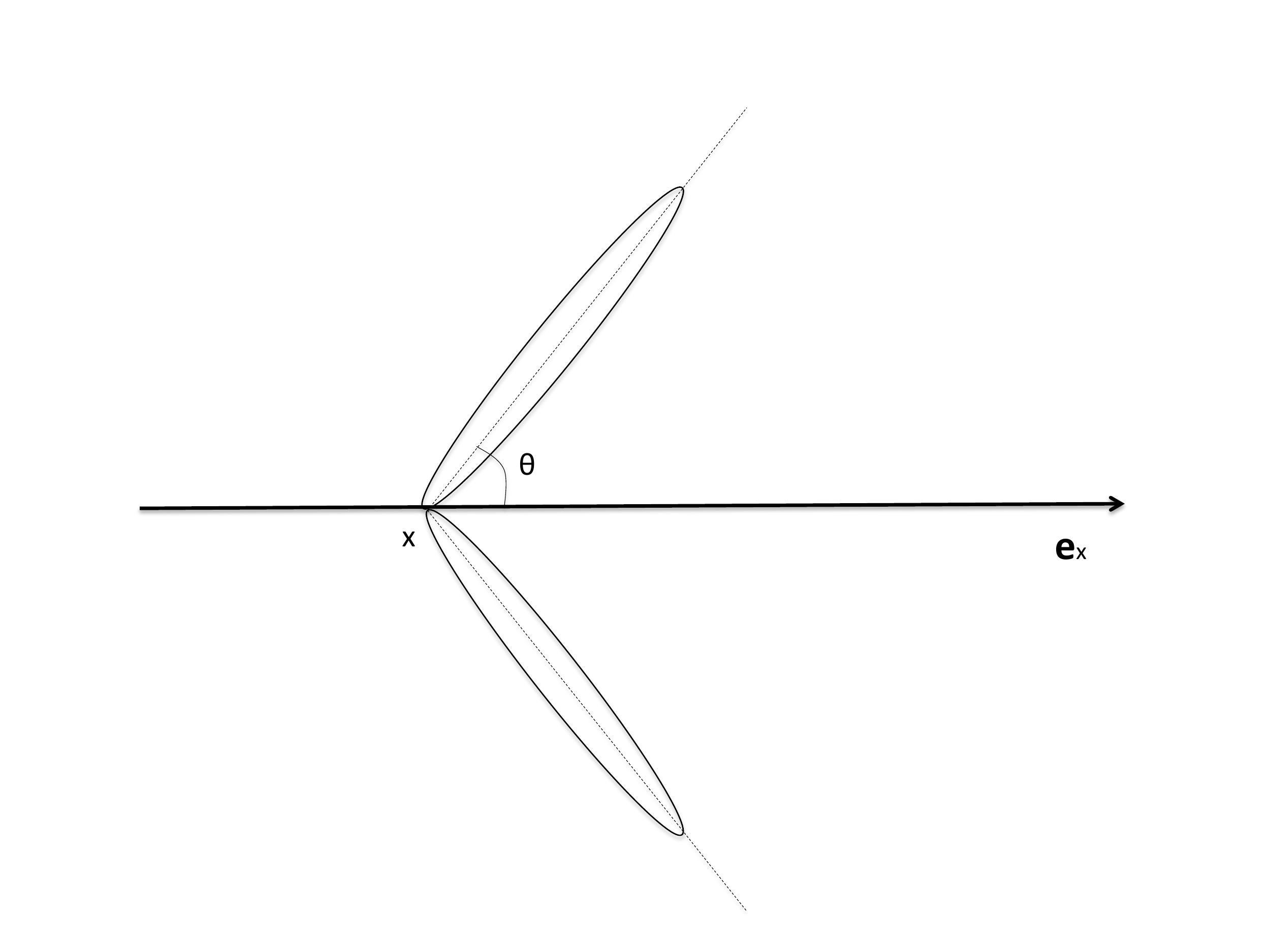}
\caption{The scallop configuration}
\label{fig:Scallop}       
\end{figure}

The possible translation of the system is determined by the consecutive opening and closing of the valves. Our aim is to determine the net translation of the body, given the function of time describing the angular velocity $\dot\theta$.

\subsection{Viscous fluid}
Here we focus on the case in which the scallop is immersed in a viscous fluid. In this regime the viscous forces dominates the inertial ones that can be neglected, so the equations governing the dynamics of the fluid are the Stokes ones:
$$
\Delta v-\nabla p=0
$$
together with the incompressibility condition $\textrm{div}\ v=0$. 
Let us consider that the ellipses have major axis $2a$ and minor axis $2b$ with $b<<a$, moreover let us suppose that $\theta\in(0,\frac{\pi}{2})$ so that it remains acute. 
One of the main difficulties in computing explicitly the equation of motion is the complexity of the hydrodynamic forces exerted by the fluid on the swimmer as a reaction to its shape changes. Since in our assumptions the minor axis of the ellipse is very small with respect to the major one, i.e. $b<<a$, we can consider the swimmer as one-dimensional, composed essentially by two links of length $2a$ (see Fig \ref{fig:Scallop}). In the case of slender swimmers, Resistive Force Theory (RFT) \cite{GrayHancock55} provides a simple and concise way to compute a local approximation of such forces, and it has been successfully used in several recent studies, see for example \cite{BeckerKoehler03,FriedrichRiedel-Kruse10}. From now on we use this approach as well, in order to obtain the forces acting on the swimmer, neglecting  the interaction between the valves. 
Since the scallop is immersed in a viscous fluid the inertial forces are negligible with respect to the viscous ones, therefore the dynamics of the swimmer follows from Newton laws in which inertia vanishes:
\begin{equation}
\vec{F}=0
\end{equation}
where $\vec{F}$ is the total force exerted on the swimmer by the fluid. As already said we want to couple the fluid and the swimmer, using the local drag approximation of Resistive Force Theory.
We denote by $s$ the arc length coordinate on the $i$-th link ($0\leq s \leq 2 a$) measured from the juncture point and
by $\vec{v}_i(s)$ the velocity of the corresponding point.
We also introduce the unit vectors $\vec{e}_1=\left(\begin{array}{cc} \cos(\theta)\\ \sin(\theta) \\ \end{array} \right)$,  $\vec{e}_1^{\bot}=\left(\begin{array}{cc} -\sin(\theta)\\ \cos(\theta) \\ \end{array} \right)$, and  $\vec{e}_2=\left(\begin{array}{cc} \cos(\theta)\\ -\sin(\theta) \\ \end{array} \right)$,  $\vec{e}_2^{\bot}=\left(\begin{array}{cc} -\sin(\theta)\\ -\cos(\theta) \\ \end{array} \right)$ in the directions parallel and perpendicular to each link and write  the position of the point at arc length $s$ as $\vec{x}_i(s)=\left(\begin{array}{c}x\\0\end{array} \right)+s\vec{e}_i$ where x is the coordinate of the joint between the two valves. 
By differentiation, we obtain,
\begin{equation}
\vec{v}_i(s) = \left(\begin{array}{c}\dot{x}\\0\end{array} \right)+ s  \dot{\theta}_i  \vec{e}_i^{\bot}\,.
\label{speed_i}
\end{equation}
The density of the force $\vec{f}_{i}$ acting on the $i$-th segment is assumed to depend linearly on the velocity. It is defined by
\begin{equation}
\vec{f}_i(s) :=-\xi \left( \vec{v}_i(s) \cdot \vec{e}_i \right) \vec{e}_i- \eta \left(\vec{v}_i(s) \cdot \vec{e}_i^{\bot}\right) \vec{e}_i^{\bot},
\label{ForceByResistiveTheory}
\end{equation}
where $\xi$ and $\eta$ are respectively the drag coefficients in the directions of $\vec{e}_i$ and $\vec{e}_i^{\bot}$ measured in $N \,s\, m^{-2}$.
We thus obtain
\begin{equation}
\vec{F}=\int_0^{2a} \vec{f}_1(s)\,ds+\int_0^{2a} \vec{f}_2(s)\,ds=0
\end{equation}
Using (\ref{speed_i}) and (\ref{ForceByResistiveTheory})  and since we are neglecting inertia we have

\begin{equation}
\label{Newton}
\begin{cases}
F_x=-4 a \xi\dot{x}\cos^2(\theta)-4 a \eta \dot{x}\sin^2(\theta)+4a^2 \eta \dot{\theta}\sin(\theta)=0\\
F_y=0
\end{cases}
\end{equation}
Observe that $F_y$ vanishes since the scallop is symmetric with respect to the $\vec{e}_x$ axis. From \eqref{Newton} is now easy to determine the evolution of $x$
\begin{equation}
\label{x_viscous}
\dot{x}=V_1(\theta)\dot{\theta}=\frac{a \eta\sin(\theta)}{\xi\cos^2(\theta)+\eta\sin^2(\theta)}\dot{\theta}
\end{equation}

\subsection{Ideal Fluid}
\label{subsec:1}
While in the previous subsection we faced the problem of the self-propulsion of the scallop immersed in a viscous fluid, here we focus on the case in which it is immersed in an ideal inviscid and irrotational fluid. Let us make the same assumptions on the parameters $a$ and $b$ that have been done in the previous section, moreover let us denote by $\Omega$ the region of the plane occupied by the swimmer in a reference configuration.\\

\noindent Assigning $(x,\theta)$ as functions of time let us call 
\begin{equation*}
\begin{aligned}
f^{(x,\theta)}:\ &\Omega\to\mathbb{R}^2\\
&\zeta\mapsto f^{(x,\theta)}(\zeta)
\end{aligned}
\end{equation*}
the function which maps each point of the swimmer $\zeta\in\Omega$ in $ f^{(x,\theta)}(\zeta)$ that is its position in the plane at time $t$. Supposing that $\theta$ can be assigned and that there are not other external forces, our aim is to find equations that describe the motion of $x$. To this end we call $v$ the velocity of the fluid, its motion is given by the Euler equations for ideal fluids
\begin{equation}
v_t+v\cdot\nabla v=-\nabla p
\end{equation}
with the incompressibility condition $\textrm{div} \ v=0$. Moreover we impose a Neumann boundary condition, that is that the normal component of the velocity of the fluid has to be equal to the normal component of the velocity of the body.
$$
\Bigl\langle v(f^{(x,\theta)})-\bigl(\frac{\partial f^{(x,\theta)}}{\partial x}\dot{x}+\frac{\partial f^{(x,\theta)}}{\partial \theta}\dot{\theta}\bigr),n^{(x,\theta)}\Bigr\rangle=0
$$
where $\Bigl\langle \cdot\Bigr\rangle$ denotes the scalar product, $n^{(x,\theta)}$ is the external normal to the set $f^{(x,\theta)}(\Omega)$. To find the evolution of $x$ we should solve the Lagrange equation
\begin{equation}
\frac{d}{dt}\frac{\partial T^b}{\partial\dot{x}}=\frac{\partial T^b}{\partial x}+F
\end{equation}

where $T^b$ is the kinetic energy of the body and $F$ the external pressure force acting on the boundary of the swimmer. As already done in \cite{Bressan07,MasonBurdick99,MunnierChambrion10} this force $F$ can be reinterpreted as a kinetic term, precisely thanks to the fact that we are in an ideal fluid. Therefore the system body + fluid is geodetic with Lagrangian given by the sum of the kinetic energy of the body ($T^b$) and the one of the fluid ($T^f$): 
$$
T^{tot}=T^b+T^f
$$
The kinetic energy of the body is the sum of the kinetic energy of the two ellipses, that reads
\begin{equation}
T^b=m\bigl(\dot{x}^2+a^2\dot{\theta}^2-2a\dot{x}\dot{\theta}\sin\theta\bigr)+I\dot{\theta}^2
\end{equation}
Since we are dealing with an ideal fluid and thus inertial forces dominates over the viscous ones, in order to derive the kinetic energy of the fluid  we will make use of the concept of \textit{added mass}. In fluid mechanics, added mass or virtual mass is the inertia added to a system because an accelerating or decelerating body must move (or deflect) some volume of surrounding fluid as it moves through it. Added mass is a common issue because the object and surrounding fluid cannot occupy the same physical space simultaneously \cite{Bessel28}. For simplicity this can be modeled as some volume of fluid moving with the object, though in reality "all" the fluid will be accelerated, to various degrees.\\
Therefore the kinetic energy of the fluid will be given by the sum of the kinetic energy of the added masses of the two ellipses.
\begin{equation}
T^f=\frac{1}{2}\vec{v}_1^T\vec{M}_{1_{add}}\vec{v}_1+\frac{1}{2}\vec{v}_2^T\vec{M}_{2_{add}}\vec{v}_2
\end{equation}
where $\vec{M}_{i_{add}}$ are the added mass matrices relative to each ellipse which are diagonal, and $\vec{v}_i$ the velocities of their centre of mass, expressed in the frame solidal to each ellipse with axes parallel and perpendicular to the major axis. 
Finally we can compute the total kinetic energy of the coupled system body+ fluid that is
\begin{equation}
\label{Tot_kin_en}
\small
\begin{aligned}
T^{tot}&=m\bigl(\dot{x}^2+a^2\dot{\theta}^2-2a\dot{x}\dot{\theta}\sin\theta\bigr)+I\dot{\theta}^2+\\
&+m_{11}\dot{x}^2\cos^2\theta+m_{22}\bigl(\dot{x}^2\sin^2\theta+a^2\dot{\theta}^2-2a\dot{x}\dot{\theta}\sin\theta\bigr)+\\
&+m_33\dot{\theta}^2
\end{aligned}
\end{equation}
Following a procedure introduced by Alberto Bressan in \cite{Bressan07}, in order to end up with a control system we perform a partial legendre transformation on the kinetic energy defining
\begin{equation*}
\begin{aligned}
p=&\frac{\partial T^{tot}}{\partial\dot{x}}=\\
&2\dot{x}\bigl(m+m_{11}\cos^2\theta+m_{22}\sin^2\theta\bigr)-2 a\dot{\theta}\sin\theta(m+m_{22})
\end{aligned}
\end{equation*}
from which we derive
\begin{equation}
\dot{x}=\frac{p+2 a\dot{\theta}\sin\theta(m+m_{22})}{2(m+m_{11}\cos^2\theta+m_{22}\sin^2\theta)}
\end{equation}
There is a wide spread literature regarding the computation of added masses of planar contours moving in an ideal unlimited fluid. We will use in the rest of the paper the added mass coefficients for the ellipse computed in \cite{Marhydro}: the added mass in the direction of the major axis is $m_{11}=\rho\pi b^2$, the one along the minor axis is $m_{22}=\rho\pi a^2$.
Notice now that writing the Hamilton equation relative to $p$, and recalling \eqref{Tot_kin_en}
$$
\dot{p}=\frac{\partial T^{tot}}{\partial x}=0
$$
thus, if we start with $p(0)=0$, $p$ remains null for all times and the evolution of $x$ becomes
\begin{equation}
\label{x_ideal}
\dot{x}=V_2(\theta)\dot{\theta}=\frac{a\sin\theta(m+\rho\pi a^2)}{m+\rho\pi b^2\cos^2\theta+\rho\pi a^2\sin^2\theta} \dot{\theta}
\end{equation}

\begin{theorem}[Scallop Theorem]
Consider a swimmer dynamics of the type
\begin{equation}
\label{linear_dyn}
\dot x= V(\theta)\dot\theta
\end{equation}
Then for every $T$-periodic deformation (i.e. stroke) one has
\begin{equation}
\Delta x=\int_0^T \dot{x}(t)\,dt=0
\end{equation}
that is, the final total translation is null
\end{theorem}
\proof{
Define the primitive of $V$ by
\begin{equation}
\label{primitive_V2}
F(\theta)=\int_0^\theta V(\sigma)\,d\sigma
\end{equation}
Then using \eqref{linear_dyn}
\begin{equation*}
\begin{aligned}
\Delta x&=\int_0^TV(\theta(t))\dot{\theta}(t)\,dt=\\
&\int_0^T\frac{d}{dt}F(\theta(t))\,dt=F(\theta(T))-F(\theta(0))=0
\end{aligned}
\end{equation*}
by the periodicity of $t\to\theta(t)$.\qed
}
\bigskip
\textit{Note that the dynamics \eqref{x_viscous} and \eqref{x_ideal} are of the type \eqref{linear_dyn}, therefore the scallop theorem is valid either in the viscous and in the ideal case.}

\section{Controllability}
\label{sec:3}
In this section we will give two different strategies to overcome the scallop theorem, both based on a switching mechanism. In particular we produce some partial and global controllability results for this switching systems.

\subsection{Partial controllability in $x$}

We have previously seen that if our scallop is immersed either in an ideal fluid or in a viscous one, if it experiences periodical shape changes it is not able to move after one cycle. Here we would like to find a way to overcome this problem. The main idea is to be able to change the dynamics during one periodical stroke and see if in this way we obtain a net motion and in particular some controllability. In order to do this we have to introduce the \textit{Reynolds number}, a number which characterizes the fluid regime. It arises from the adimesionalization of the Navier-Stokes equations and it is defined by
\begin{equation}
Re=\frac{VL\rho}{\eta}=\frac{VL}{\nu}
\end{equation}
 where $V$ is the characteristic velocity of the body immersed in the fluid, $L$ its characteristic length, $\rho$ the density of the fluid, $\eta$ its viscosity and  $\nu=\frac{\eta}{\rho}$ is the kinematic viscosity. The Reynolds number quantifies the relative importance of inertial versus viscous effects. 
 
 \subsubsection{$\boxed{\eta=\eta(|\dot\theta|)}$}
 
\begin{figure}[H]
\begin{center}
\includegraphics[width=%
0.4\textwidth]{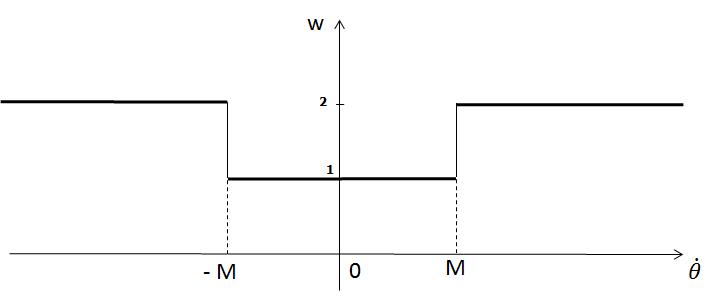} 
\end{center}
\caption{The rule of the classical switching}
\label{Fig1}
\end{figure}

 Let us recall that if $v(t,x)$ is a solution of the Navier Stokes equations, the function $u(t,x)=v(ct,x)$, $c>0$ is still a solution of the Navier Stokes equations but with a different viscosity.
 Now assume that the absolute value of the speed $\dot{\theta}$ is very high, this means that rescaling the time of the solution of the Navier Stokes equations, we end up with a viscosity $\eta$ that is very small and therefore the Reynolds number is large. In this case the inertial forces dominates over the viscous ones, so we can consider the scallop immersed in an ideal fluid and thus use the dynamics \eqref{x_ideal}. Then we suppose that at a certain point of the cycle the absolute value of the angular velocity is very small. In this case we have a solution of the Navier Stokes equations with a very high viscosity $\eta$. Thus we can suppose that the scallop is immersed in a Stokes fluid, since the viscous effects dominates the inertial ones and use the dynamics \eqref{x_viscous}. This situation is well represented by a switching system in which the change of the dynamics is determined by the modulus of the angular velocity $\dot{\theta}$: if it is big (i.e $|\dot\theta|>M$ with $M>0$) we use the ideal approximation and the corresponding dynamics; if it is small (i.e $|\dot\theta|<M$ with $M>0$) we use instead the viscous approximation and the relative dynamics.The switching rule in Fig \ref{Fig1} should also consider what happens when $|\dot\theta|=M$. However in the sequel we are going to exhibit a function $\dot\theta$ which stays in $M$ or $-M$ for only a set of times of null measure.\\
Our aim is to prove that using this kind of switching we are able to have a net displacement, both forward or backward, using periodic continuous functions $\dot{\theta}$

 According to what said before we can prescribe the angular velocity $\dot\theta$ and thus use it as a control function $u$. Therefore we write the system as a control system that is
\[
\begin{cases}
\dot{x}(t) = V_{w(t)}(\theta(t))u(t), \\
\dot \theta(t)=u(t)\\
w(t)= h[u](t) \\
x(0)=x_0, \ \theta(0)=\theta_0\ \ w(0)= w_0
\end{cases} \]
where $u$ is continuous and 
$$
h[u]=
\begin{cases}
2\quad\text{if $|u|>M$}\\
1\quad\text{if $|u|<M$}
\end{cases}
$$

Moreover let us call  $F_i$ the primitives of the functions $V_i$, for $i=1,2$. They are :
\begin{equation*}
F_1 =\frac{a \eta\arctan h(\sqrt{\frac{\eta-\xi}{\eta}}\cos\theta)}{\sqrt{\eta(\eta-\xi)}},
\end{equation*}
\begin{equation*}
F_2 = \frac{-a\sqrt{m+a^2\rho\pi}\arctan h(\frac{\sqrt{(a^2-b^2)\rho\pi}\cos\theta}
{\sqrt{m+a^2\rho\pi}})}{\sqrt{\rho\pi(a^2-b^2)}}
\end{equation*}

\begin{theorem}
\label{modulus}
With the previous switching scheme we are able to overcome the Scallop paradox, thus to move both forward and backward. More precisely there are $r>0$ small enough (see remark \ref{erre1}),  a final time $T>0$ and a continuous $T$-periodic control function $u(t)$, which make the system move between two fixed points along the $x$ axis, $x_0$ and $x_f\in]x_0-r,x_0+r[$, in the time $T$.
\end{theorem}
\proof{
\textbf{First case: $u(0)>M$} \\In this case we start with the ideal approximation (i.e $w_0=2$)
\begin{equation}\label{caso1}
V_{w(t)}(\theta(t))=
\begin{cases}
V_2(\theta(t)) & 0<t<t_1,\\
V_1(\theta(t)) & t_1<t<t_2 \\
V_2(\theta(t)) & t_2<t<t_3\\
V_1(\theta(t)) & t_3<t<t_4 \\
V_2(\theta(t)) & t_4<t<T\\
\end{cases}
\end{equation}
 with
 \begin{align*}
& t_1:=inf\{T>t>0\,|\, u(t)=M\}\quad \text{ and } \\
  &t_2:=inf\{T>t>t_1\,|\,u(t)=-M\} \quad \text{ and } \\
  & t_3:=inf\{T>t>t_2\,|\, u(t)=-M\}\quad \text{ and } \\
  & t_4:=inf\{T>t>t_3\,|\, u(t)=M\}\
 \end{align*}
 assuming that $inf \, (\emptyset)=+\infty$. 
The net motion is then calculated as
\begin{equation}\label{Spost0caso3}
\begin{aligned}
\Delta x= &(F_2-F_1)(\theta(t_1))+(F_2-F_1)(\theta(t_3))\\
&-(F_2-F_1)(\theta(t_2))-(F_2-F_1)(\theta(t_4)).
\end{aligned}
\end{equation}
taking into account that $\theta(0)=\theta(T)$ and that $(F_2-F_1)(\theta(t_i))$ does not appear in the equation if $t_i=+\infty$.\\
We want to prove that we are able to move  choosing a suitable periodic evolution for our control function $\dot \theta=u$. 
Let us call the unknowns $\theta_i:=\theta(t_i)$, for $i=1\dots 4$. First of all we show that $\Delta x$ as function of $(\theta_1,\theta_2,\theta_3,\theta_4)$ is surjective in $]0,\frac{\pi}{2}[\times]0,\frac{\pi}{2}[\times]0,\frac{\pi}{2}[\times]0,\frac{\pi}{2}[$.\\
We are going to prove that
\begin{equation*}
\begin{aligned}
&\nabla(\Delta x)=\left(\begin{array}{cc}- (V_2-V_1)(\theta_1)\\ (V_2-V_1)(\theta_2) \\- (V_2-V_1)(\theta_3) \\ (V_2-V_1)(\theta_4)\end{array} \right)\neq 0 \\
& \text{in} \ \ (\theta_1,\theta_2,\theta_3,\theta_4) \in ]0,\frac{\pi}{2}[\times]0,\frac{\pi}{2}[\times]0,\frac{\pi}{2}[\times]0,\frac{\pi}{2}[
\end{aligned}
\end{equation*}
so that \eqref{Spost0caso3} is a submersion and surjective as required.\\ 
Recall that the function $(F_2-F_1)(\cdot)$ is always increasing indeed
\begin{equation}
\small
\begin{aligned}
&\frac{(F_2-F_1)(\theta)}{\partial\theta} = \\
&\bigg( -\frac{a\eta}{\xi\cos^2\theta + \eta\sin^2\theta}+\frac{ma+\rho\pi a^2}{m+\rho\pi b^2 cos^2\theta+\rho\pi a^2\sin^2\theta}\bigg)\sin\theta\\
&\\
&=\frac{\sin\theta \cos^{2}\theta\big( ma(\eta-\xi)+\rho\pi(\xi a^2-\eta b^2)\big)}{(m+\rho\pi b^2 cos^2\theta+\rho\pi a^2\sin^2\theta)(\xi\cos^2\theta + \eta\sin^2\theta)}> 0 \\\\&\qquad\qquad \text{for}\ \ \theta \in ]0, \frac{\pi}{2}[\quad\text{and }b<<a
\end{aligned}
\end{equation}
From this immediatly follows that $\nabla(\Delta x)\neq 0$.\\
The surjectivity ensures us that for any fixed $\Delta x$ in a neighborhood of zero we are always able to find a $(\theta_1,\theta_2,\theta_3,\theta_4)$ which realize the desired displacement. Moreover, thanks to the symmetry properties of the function defining the displacement, also each of the $4$-uplets $(\theta_1,\theta_4,\theta_3,\theta_2)$, $(\theta_3,\theta_4,\theta_1,\theta_2)$   and $(\theta_3,\theta_2,\theta_1,\theta_4)$ realizes the same displacement. 
Supposing $\Delta x>0$ and recalling that the function $(F_2-F_1)(\cdot)$ is increasing, then the angles $(\theta_1,\theta_2,\theta_3,\theta_4)$ will have a suitable order that can or not be coherent with the switching rule and the periodicity of $\dot\theta$. If their sorting is appropriate we will choose a control $\dot\theta=u$ such that $\theta(t_i)=\theta_i$. Otherwise at least one of the $4$-uplets above will be right. Thus defining $(\theta^{'}_1,\theta^{'}_2,\theta^{'}_3,\theta^{'}_4)$ this latter uple, we take a control $u$ such that $\theta(t_i)=\theta^{'}_i$. This choice of the control will lead us to obtain the desired positive displacement.\\

For example suppose that the uplet $(\theta_1,\theta_2,\theta_3,\theta_4)$ which realizes the desired positive displacement, satisfy $\theta_3>\theta_4>\theta_1>\theta_2$. Indeed
\begin{align*}
&(F_2-F_1)(\theta_1)-(F_2-F_1)(\theta_2)>0\\
&&\Longrightarrow\Delta x>0\\
&(F_2-F_1)(\theta_4)-(F_2-F_1)(\theta_3)<0
\end{align*}
To respect the switching scheme in the time interval $(t_2,t_3)$ the function $\dot\theta$ should decrease and thus $\theta_2>\theta_3$. The latter is not satisfied by $(\theta_1,\theta_2,\theta_3,\theta_4)$, but taking $(\theta^{'}_1,\theta^{'}_2,\theta^{'}_3,\theta^{'}_4)=(\theta_3,\theta_4,\theta_1,\theta_2)$, we have the same $\Delta x$ and the switching scheme is now respected. 
\begin{figure}[H]
\begin{center}
\includegraphics[width=%
0.55\textwidth]{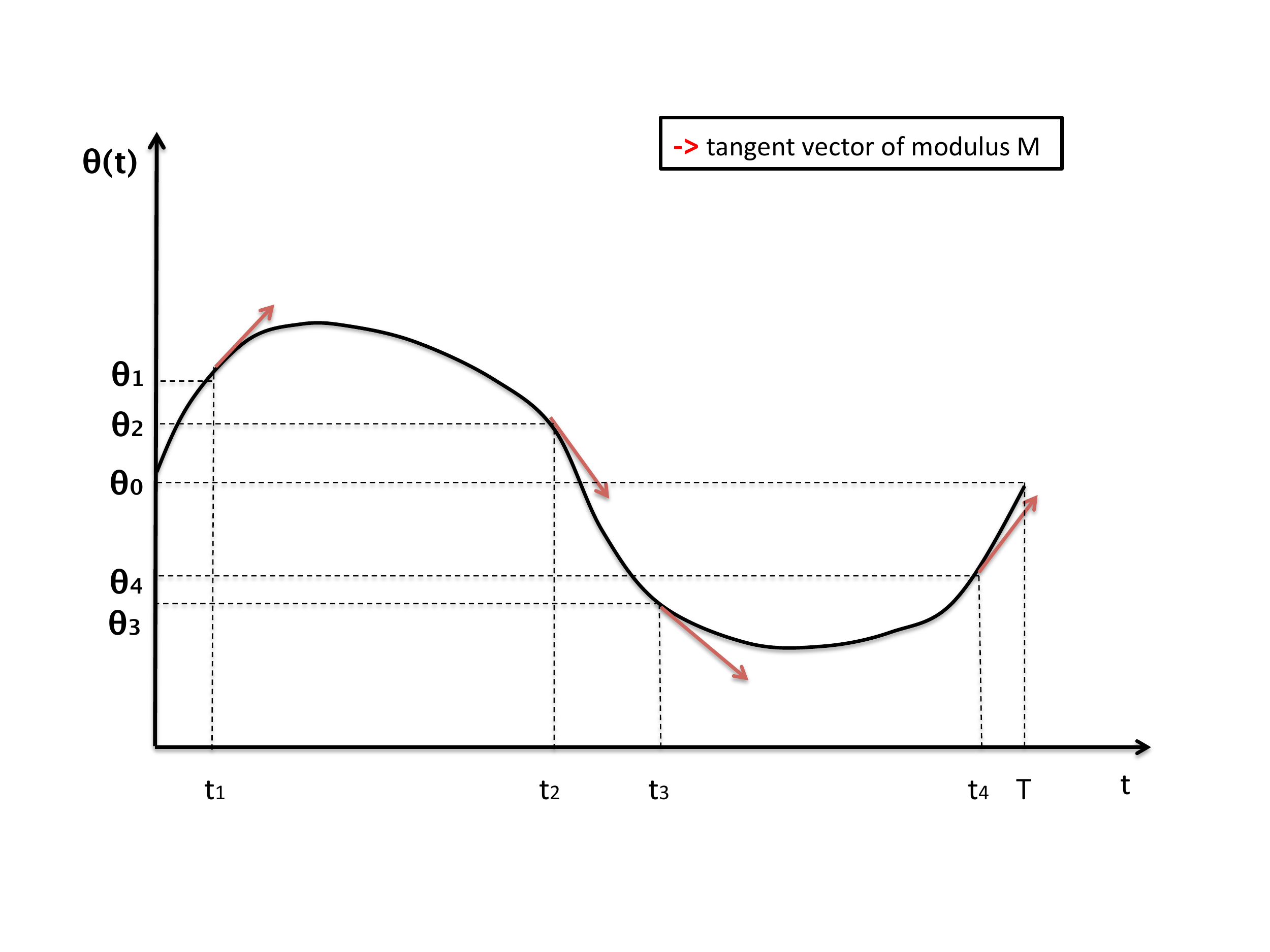} 
\end{center}
\caption{This figure shows a possible choice of $\theta(t)$ which realizes a positive displacement and respects the switching scheeme}
\label{Theta_rele3}
\end{figure}
Analogous arguments can be used if $\Delta x<0$.

\textbf{Second case: $-M<u(0)<M$}\\ In this case we start form the viscous approximation  (i.e $w_0=1$).
Using arguments similar to the ones used before to compute $\Delta x$ and to prove its surjectivity, and redefining accordingly the times $t_i$ for $i=1\dots 4$. we have that
\begin{equation}\label{Spost0caso3_2}
\begin{aligned}
\Delta x= &(F_2-F_1)(\theta(t_2))+(F_2-F_1)(\theta(t_4))\\
&-(F_2-F_1)(\theta(t_1))-(F_2-F_1)(\theta(t_3)).
\end{aligned}
\end{equation}
going on as before, exploiting the surjetivity and the symmetry of the last function, we are able to find a control $u$ that realizes the desired displacement.

\textbf{Third case: $u(0)<-M$}\\ This case is analogous to the first one.\\
In  conclusion we have proved that wherever we start on the switching diagram we are able to achieve  a net displacement either positive or negative and then we have the controllability.\qed
}

 \begin{remark}
 \label{erre1}
 Note that the value of $r$ in the the last theorem is the maximal value that the function\\ $|\Delta x(\theta_1,\theta_2,\theta_3,\theta_4)|$ can assume in \\$]0,\frac{\pi}{2}[\times]0,\frac{\pi}{2}[\times]0,\frac{\pi}{2}[\times]0,\frac{\pi}{2}[$. Thus the constant $r$ is independent from $x$ and $\theta$.\\ To cover distances $|\Delta x|\geq r$ we should divide the spatial interval in $N$ subintervals of length less than $r$, each one realized by a $u$ of period $\frac{T}{N}$. Repeating $N$ times this control $u$ we are able to reach the desired displacement. 
 \end{remark}

\subsubsection{$\boxed{\eta=\eta(sign(\dot\theta))}$}

\begin{figure}[H]
\begin{center}
\includegraphics[width=%
0.4\textwidth]{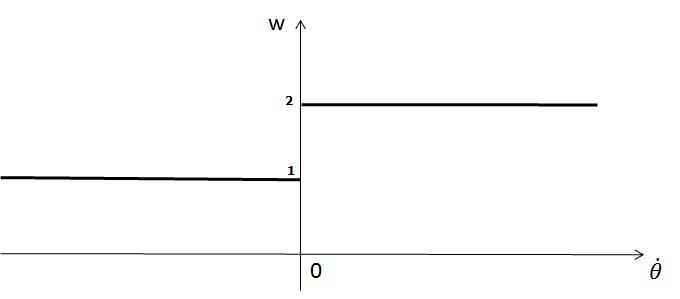} 
\end{center}
\caption{The rule of the classical switching}
\label{classical}
\end{figure}

While in the previous subsection we supposed that the change in the fluid regime was linked to the magnitude of the modulus of the angular velocity, here we would like to link the two fluids approximations to the sign of $\dot\theta$. This model can be used to describe a different response of the fluid to the opening and closure of the scallop' valves. For instance  the situation in which the fluid has a pseudoelastic nature that assists the valve opening but resist the valve closing \cite{Pseudoelastic}. Thus, according to our assumption the viscosity of the fluid changes between the opening and the closing of the valves, switching from one constant value to another one.\\
This can be represented by a switching scheme as in Fig. \ref{classical}. If the valves are opening ($\dot\theta>0$) we suppose that the fluid is not opposing resistance, as if the scallop is immersed in an ideal fluid; instead when the valves are closing  ($\dot\theta<0$) the fluid is opposing a big resistance, and we can consider the scallop immersed in a viscous fluid. 

The system can be written as a control system, in which the control function $u(t)$ is the angular velocity $\dot\theta$:
\[
\begin{cases}
\dot{x}(t) = V_{w(t)}(\theta(t))u(t), \\
\dot \theta(t)=u(t)\\
w(t)= h[u](t) \\
x(0)=x_0, \ \theta(0)=\theta_0\ \ w(0)= w_0
\end{cases} \]
where the control $u$ is continuous and now
$$
h[u]=
\begin{cases}
2\quad\text{if $u>0$}\\
1\quad\text{if $u<0$}
\end{cases}
$$

\begin{theorem}
With the classical switching scheme (see Fig \ref{classical}) we are able to overcome the scallop theorem but moving only forward. That is, there are $r>0$ small enough, a time $T>0$ and a continuous $T$-periodic control function, which make the system move between two fixed configurations $x_0$ and $x_f$ with $x_f\in [x_0, x_0+r[$, in the time $T$.
\end{theorem}

\begin{figure}[H]
\begin{center}
\includegraphics[width=%
0.4\textwidth]{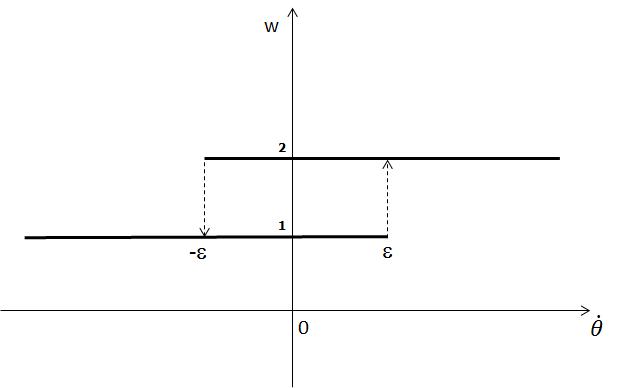} 
\end{center}
\caption{The thermostatic approximation}
\label{Fig00}
\end{figure}

\proof{
Let us suppose to start with the ideal approximation, so that we are opening the valves
\begin{equation*}
u(0)> 0 \ \ \ \text{and} \ \ \ w_0=2,
\end{equation*}
\begin{equation}\label{caso3}
V_{w(t)}(\theta(t))=
\begin{cases}
V_2(\theta(t)) & 0<t<t_1,\\
V_1(\theta(t)) & t_1<t<t_2, \\
V_2(\theta(t)) & t_2<t<T.
\end{cases}
\end{equation}
with 

 \begin{align*}
& t_1:=inf\{T>t>0\,|\, u(t)=0\}\quad \text{ and } \\
  &t_2:=inf\{T>t>t_1\,|\,u(t)=0\}
 \end{align*}
 with $inf(\emptyset)=+\infty$.
The net motion can be computed as
\begin{equation}\label{Spost0caso}
\Delta x= F_2(\theta(t_1))+F_1(\theta(t_2))-F_1(\theta(t_1))-F_2(\theta(t_2)).
\end{equation}
recalling as before that $\theta(0)=\theta(T)$.
We want to prove that we are able to move  choosing a suitable periodic evolution for our control function $\dot \theta=u$. 
Let us call $\theta_1:=\theta(t_1)$ and $\theta_2:=\theta(t_2)$, first of all we show that $\Delta x$ as function of $(\theta_1,\theta_2)$ is surjective in $]0,\frac{\pi}{2}[\times]0,\frac{\pi}{2}[$.\\
Like before we prove that 
\begin{equation*}
\begin{aligned}
&\nabla(\Delta x)=\left(\begin{array}{cc} (V_2-V_1)(\theta_1)\\ (V_1-V_2)(\theta_2) \\ \end{array} \right)\neq 0 \\
& \text{in} \ \ (\theta_1, \theta_2) \in ]0, \frac{\pi}{2}[\times]0, \frac{\pi}{2}[
\end{aligned}
\end{equation*}
hence \eqref{Spost0caso} is a submersion and surjective as required. 
Notice that 
$$
\Delta x=(F_2-F_1)(\theta_1)-(F_2-F_1)(\theta_2)
$$
If we chose a control such that $\theta_1>\theta_2$ then $\Delta x$ will be positive, while if $\theta_1<\theta_2$ then $\Delta x$ will be negative. But since we need to respect the switching rule the last case could not be achieved because after $t_1$ $\dot\theta=u<0$ and thus we are closing the valves therefore $\theta(t_2)=\theta_2$ will be necessarily less than $\theta(t_1)=\theta_1$.\\
 The case where $u(0)<0$ is analogous to the previous one.\\
 In conclusion we have proved that for every choice of $w_0$ we are able to achieve a net displacement but only forward. \qed}
\paragraph{\textit{Thermostatic case\\}}
In this section we introduce a variant of the previous switching in order to be able to move both forward and backward and therefore have a result of partial controllability in $x$. Our approach is to link the variation of $u=\dot{\theta} \in \mathbb{R}$ by a delayed thermostat, an operator with memory, introduced rigorously in \cite{Visintin}, consisting of two different thresholds for passing separately from one edge to the another one and vice-versa. This idea was inspired by \cite{Nature} in which the Scallop opening and closing is actuated by an external magnetic field, and thus a delay mechanism is reasonable. We suppose that the dynamics $V$  depends on the angle  $\theta \in ]0, \frac{\pi}{2}[$, and also depends on a discrete variable $w \in \left\{1, 2\right\}$, whose evolution is governed by a delayed thermostatic rule, subject to the evolution of the control $u$. In Fig. \ref{Fig00} the behavior of such a rule  is explained, correspondingly to the choice of a fixed threshold parameter $\varepsilon>0$. The output $w\in\{1,2\}$ may jump  from $2$ to $1$ only when the input $u$ is equal to $-\varepsilon$, and must jump when $U$ coming from the right (i.e. from values larger than or equal to $-\varepsilon$), possibly goes below the threshold $-\varepsilon$; it may jump from $1$ to $2$ only when $u$ is equal to $\varepsilon$, and must jump when it comes from the left (i.e. from values smaller than or equal to $\varepsilon$) possibly goes above the threshold $\varepsilon$. In all other situations it remains locally constant in time. In particular, when $u>\varepsilon$ then $w$ is equal to $2$, and when $u<-\varepsilon$ then $w$ is equal to $1$.

The controlled evolution is then given by
\begin{equation}
\label{control_syst_scallop_hyst}
\begin{cases}
\dot{x}(t) = V_{w(t)}(\theta(t))u(t), \\
\dot\theta(t)=u(t)\\
w(t)= h_{\varepsilon}[u](t) \\
x(0)=x_0, \ \ \theta(0)=\theta_0 \ \ w(0)= w_0
\end{cases} 
\end{equation}
where $h_{\varepsilon}\left[\cdot\right]$ represents the thermostatic delayed relationship between the input $u$ and the output $w$. Note that the initial value $w_0 \in \left\{1, 2\right\}$ must be coherent with the thermostatic relation: $w_0 = 2$ (resp. $w_0 = 1$) whenever $\dot{\theta}_0 > \varepsilon$ (resp. $\dot{\theta}_0 < -\varepsilon$).\\
We start now to analyse the value of the displacement $\Delta x$ depending of the value of $u$ proving the following result:
\begin{theorem}
\label{Controll_hyst_x}
Let $x_f\in]x_0-r,x_0+r[$ with $r>0$ small enough. Then, there always exits a time $T>0$ and a continuous $T$-periodic control function $\dot\theta=u$ (hence a periodic $\theta$)  such that one can move from $x_0$ to $x_f$ in time $T$ 
when the delayed thermostat is taken into account. In other words the system \eqref{control_syst_scallop_hyst} is partially controllable in $x$.
\end{theorem}
\proof{
{\bf First case}
\begin{equation*}
-\varepsilon< u(0)< \varepsilon \ \ \ \text{and} \ \ \ \ w_0=1
\end{equation*}
then we have
\begin{equation}\label{Internosotto}
V_{w(t)}(\theta(t))=
\begin{cases}
V_1(\theta(t)) & 0<t<t_1\\
V_2(\theta(t)) & t_1<t<T.
\end{cases}
\end{equation}
where $t_1$ is the first time for which $u$ goes through $\varepsilon$, i.e.
$$
t_1:=inf\{T>t>0\,|\, u(t)=\varepsilon\}
$$
 and $T$ is the final time. The displacement is then
\begin{equation}\label{Spost1caso}
\Delta x= F_1(\theta(t_1))-F_1(\theta(0))+F_2(\theta(0))-F_2(\theta(t_1)).
\end{equation}
recallin as before that $\theta(0)=\theta(T)$.

We call $\theta(t_1)=\theta_1$ 
and we want to prove that we able to obtain $\Delta x= c, \forall \ |c| <r$ using a suitable periodic control function. In order to do this we show that $\Delta x(\theta_1)$ is surjective in a neighborhood of zero. First of all we compute the derivative and show that it is different from $0$ and negative.
\begin{equation*}
\small
\begin{aligned}
&\frac{\partial\Delta x}{\partial\theta_1} = V_1(\theta_1)-V_2(\theta_1) =\\
&\bigg(     \frac{a\eta \sin\theta_1}{\xi\cos^2\theta_1 + \eta\sin^2\theta_1}-\frac{(ma+\rho\pi a^2)\sin\theta_1}{m+\rho\pi b^2 cos^2\theta_1+\rho\pi a^2\sin^2\theta_1}\bigg)\\
&\\
&=\frac{ \sin\theta_1 \cos^{2}\theta_1\big(-ma(\eta-\xi)-\rho\pi(\xi a^2-\eta b^2)\big)}{(m+\rho\pi b^2 cos^2\theta_1+\rho\pi a^2\sin^2\theta_1)(\xi\cos^2\theta_1 + \eta\sin^2\theta_1)}\neq 0 \\
& \text{for}\ \ \theta_1 \in ]0, \frac{\pi}{2}[
\end{aligned}
\end{equation*}
Notice also that since in our assumptions $b$ is negligible with respect to $a$, i.e $b<<a$, we have that $\eta b^2<< \xi a^2$ and thus the derivative is always negative and consequently the $\Delta x$ is decreesing.
We are interested in $ \theta_1 \in ]0, \frac{\pi}{2}[$. 
Since the derivative of the function defining the displacement is different from $0$ in $]0,\frac{\pi}{2}[$, \eqref{Spost1caso} is locally invertible. Thus, since  the inverse image of $0$ is $\theta_0$  then the inverse image of a neighborhood of $0$ is a neighborhood of $\theta_0$. Finally, recalling that $\Delta x$ is decreasing, we can conclude that \eqref{Spost1caso} can be positive or negative i.e. if we chose a control such that $\theta_1<\theta_0$ the displacement will be positive instead if $\theta_1>\theta_0$ it will be negative. In both cases the switching rule is respected thanks to the presence of the thermostat.

\begin{figure}[H]
\begin{center}
\includegraphics[width=%
0.5\textwidth]{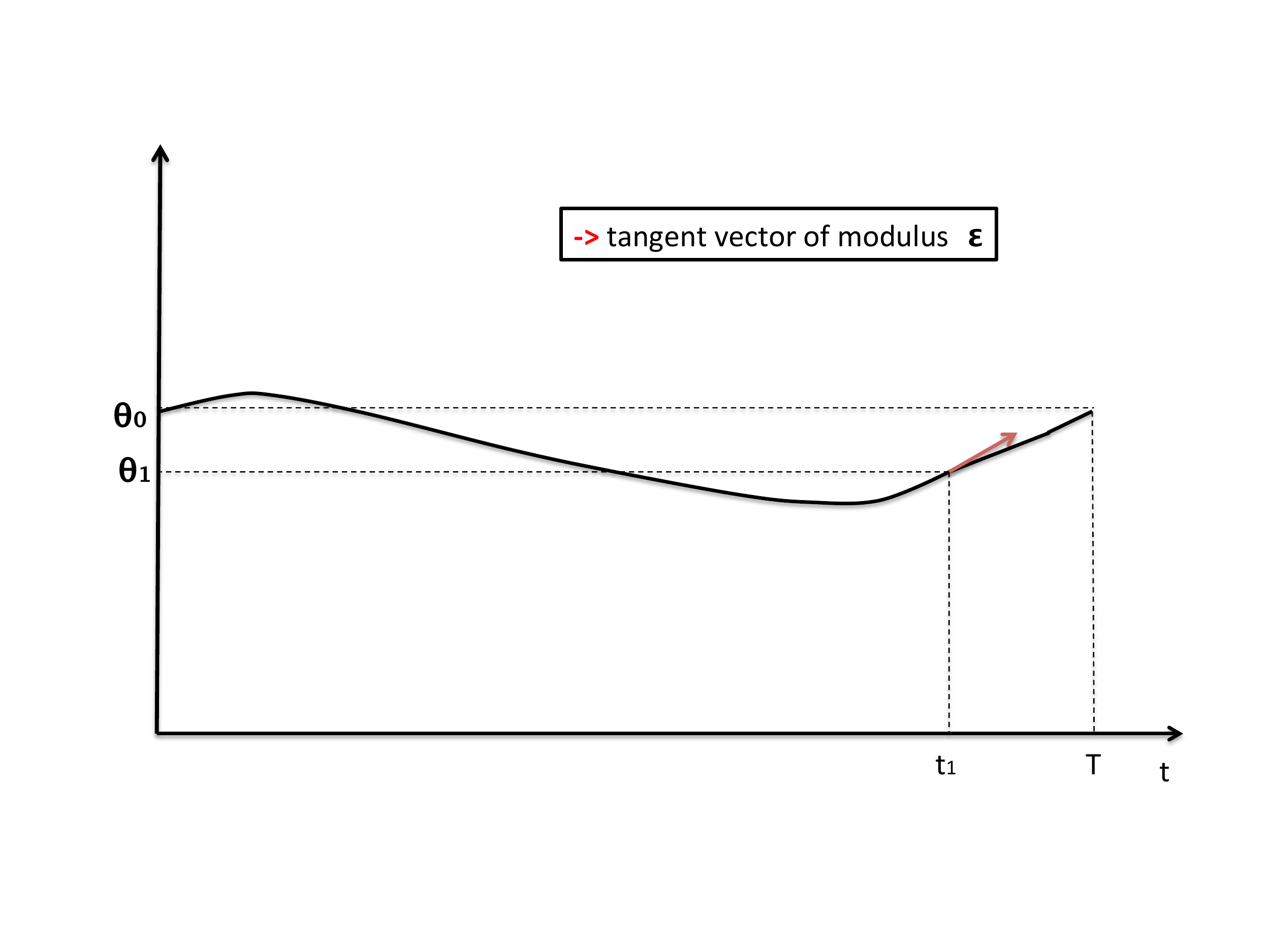} 
\end{center}
\caption{A possible choice of $\theta(t)$ starting from $0<u(0)<\varepsilon$ which realizes a positive displacement}
\label{Theta_thermostat}
\end{figure} 

{\bf Second case}
\begin{equation*}
-\varepsilon < u(0)< \varepsilon \ \ \ \text{and} \ \ \ w_0=2
\end{equation*}
then we have
\begin{equation}\label{Internosopra}
V_{w(t)}(\theta(t))=
\begin{cases}
V_2(\theta(t)) & 0<t<t_1\\
V_1(\theta(t)) & t_1<t<T.
\end{cases}
\end{equation}
where $t_1$ is the first time for which $u$ goes through $-\varepsilon$ 
$$
t_1:=inf\{T>t>0\,|\, u(t)=-\varepsilon\}
$$
and $T$ the final time. The displacement is
\begin{equation}\label{Spost2caso}
\Delta x= F_2(\theta(t_1))-F_2(\theta(0))+F_1(\theta(T))-F_1(\theta(t_1)).
\end{equation}
Calling again $\theta(t_1):=\theta_1$ also in the case we verify the surjectivity showing that the derivative of the displacement is different from zero.

Hence \eqref{Spost2caso} is locally invertible and the inverse image of a neighborhood of $0$ is a neighborhood of $\theta_0$. We can conclude as in the previous case that \eqref{Spost2caso} can be either positive or negative choosing a suitable control.\\

{\bf Third case}
\begin{equation*}
u(0)> \varepsilon \ \ \ \text{and} \ \ \ w_0=2,
\end{equation*}
\begin{equation}\label{caso3}
V_{w(t)}(\theta(t))=
\begin{cases}
V_2(\theta(t)) & 0<t<t_1,\\
V_1(\theta(t)) & t_1<t<t_2, \\
V_2(\theta(t)) & t_2<t<T.
\end{cases}
\end{equation}
with 

 \begin{align*}
& t_1:=inf\{T>t>0\,|\, u(t)=-\varepsilon\}\quad \text{ and } \\
& t_2:=inf\{T>t>t_1\,|\, u(t)=\varepsilon\}
 \end{align*}
The net motion is
\begin{equation}\label{Spost3caso}
\Delta x= F_2(\theta(t_1))+F_1(\theta(t_2))-F_1(\theta(t_1))-F_2(\theta(t_2)).
\end{equation}
recalling that $\theta(0)=\theta(T)$.\\
Also in this case we want to prove that we are able to move both forward or backward. Therefore we show that $\Delta x$ is surjective in $]0,\frac{\pi}{2}[\times]0,\frac{\pi}{2}[$ as in the non hysteretic case.
We compute the gradient and show that it is never null
\begin{equation*}
\nabla(\Delta x)=\left(\begin{array}{cc} (V_2-V_1)(\theta_1)\\ (V_1-V_2)(\theta_2) \\ \end{array} \right)\neq 0 \ \ \text{in} \ \ (\theta_1, \theta_2) \in ]0, \frac{\pi}{2}[\times]0, \frac{\pi}{2}[
\end{equation*}
hence \eqref{Spost3caso} is a submersion and surjective as required. 
Notice that 
$$
\Delta x=(F_2-F1)(\theta_1)-(F_2-F_1)(\theta_2)
$$
and recall that the function $(F_2-F1)(\cdot)$ is always increasing. Hence, if we use a control such that $\theta_1>\theta_2$ the $\Delta x$ will be positive, while if $\theta_1<\theta_2$ then $\Delta x$ will be negative. Also in this case both the alternatives can be achieved respecting the switching rule.
Therefore we are able to obtain the desired displacement.\\\\
{\bf Fourth case}
 The case where $u(0)<- \varepsilon$ is analogous to the previous one.\\\\
In conclusion we have proved that for every choice of $w_0$ we are always able to find a periodic and continuous control $\dot{\theta}=u$ that allows us to obtain the desired displacement. The system \eqref{control_syst_scallop_hyst} is then partially controllable in $x$.\qed}
 \begin{remark}
The introduction of the thermostat is essential because allows us to achieve displacements of every sign and thus the controllability result in $x$. This fact is strictly linked to the presence of the thresholds, indeed we are allowed to move between them without changing dynamics and therefore obtain values $\theta_1<\theta_2$ either $\theta_1>\theta_2$, and thus move both forward and backward.
 \end{remark}
 \begin{remark}
 \label{erre}
 Note that the maximal value of $r$ in the last theorem is $|\Delta x(\frac{\pi}{2})|$ if $-\varepsilon<u(0)<\varepsilon$, and  $|\Delta x(\frac{\pi}{2},0)|$ if $-\varepsilon<u(0)$ or $u(0)>\varepsilon$. Thus it is always independent from $x$ and $\theta$.\\ To cover distances $|\Delta x|\geq r$ we should divide the spatial interval in $N$ subintervals of length less than $r$, each one realized by a $u$ of period $\frac{T}{N}$. Repeating $N$ times this control $u$ we are able to reach the desired displacement.
 \end{remark}
\subsection{Global controllability result}

In this subsection we are interested in studying whether it is feasible for the system of the scallop to move between two fixed configurations ($(x_0,\theta_0)$ and $(x_f,\theta_f)$). This part add something to the previous one, since we are prescribing both the initial and final positions and angles. The following holds:
\begin{theorem}
\label{control_x_theta}
Let $A$ and $B$ be two fixed positions along the $x$-axis and $\theta_0$, $\theta_f$ two fixed angles. Then, we are always able to find a suitable control function $u(t)$ such that the scallop system moves between $A$ and $B$ passing from $\theta(0)=\theta_0$ to $\theta(T)=\theta_f$, where $T$ is a suitable big enough final time. Moreover such function $u(t)$ respects the switching rules modeling the dependence of  the viscosity $\eta$ from $|\dot\theta|$ (Fig.\ref{Fig1}) and from $sign(\dot\theta)$ with the thermostat, (Fig.\ref{Fig00}). In other words the system \eqref{control_syst_scallop_hyst} is controllable.
\end{theorem}
\proof{
Let $u(t)$ the periodic function that makes the system move between $A$ and $B$ with final angle $\theta_0$ during a time $t^{'}$. We have proved the existence of such a function with both switching rules, in the previous subsection. 
Now whatever $w(t^{'}) $ we open or close the valves respecting the switching rule in Fig.\ref{Fig1} or Fig. \ref{Fig00} respectively until we reach the desired angle $\theta_f$.
%
We call $t^{''}$ the time in which we have $\theta_f$ and $C$ the point in which we are arrived.
Now starting from $C$ with $w(t^{''})$  we move to $B$ using another periodic $u(t)$ (hence $\theta(t)$ periodic), whose existence is ensured from Theorem \ref{Controll_hyst_x}\qed
}
\begin{figure}[H]
\begin{center}
\includegraphics[width=%
0.4\textwidth]{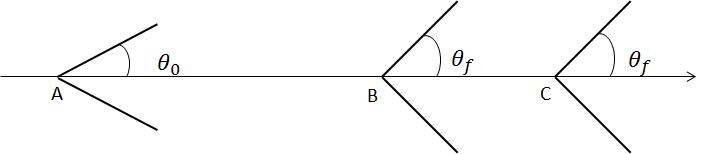} 
\end{center}
\caption{This figure represents one of the cases considered in the proof of Theorem \ref{control_x_theta} }
\label{Fig2}
\end{figure}

\section{Numerical examples}
\label{sec:4}
In this section we will show, through numerical simulations, that our theoretical pretictions on the controllability of the Scallop along $x$ are good. Moreover we will also describe how it is possible to obtain the same results removing the continuity hypothesis on $\dot\theta$.
In what follows the pictures are all relative to the controllability result which follows the thermostatic switching scheme (see Fig \ref{Fig00}) that is the most interesting one. Similar results can be obtained analogously using the other switching described in Fig. \ref{Fig1}.\\

Let us suppose to start with $w(0)=2$ which means $\dot\theta(0)>\varepsilon$, the following pictures show a possible choice of the control $\dot\theta$ to obtain a displacement $\Delta x=1\, cm$, using the following parameters: $a=2 \,cm$, $b=0.1\,cm$, $\eta=2\, N s m^{-2}$, $\xi=1\, N s m^{-2}$ $m=1\, g$ and $\rho=1 \,g cm^{-3}$. More precisely in these simulations we decided to use a periodic polynomial control $\theta(t)$ that can be uniquely determined imposing the following constraints.
\begin{equation}
\label{constraints}
\begin{aligned}
&\dot\theta(0)=\dot\theta_0&&\dot\theta(t_1)=-\varepsilon&&\dot\theta(t_2)=\varepsilon&&\dot\theta(T)=\dot\theta_0\\
&\theta(0)=\theta_0&&\theta(t_1)=\theta_1&&\theta(t_2)=\theta_2&&\theta(T)=\theta_0
\end{aligned}
\end{equation}
where $\theta_1$ and $\theta_2$ are determined by the numerical inversion of the function $\Delta x$ \eqref{Spost1caso} and we chose $t_1=2 \,s$, $t_2=6 \,s$ and $T=7 \,s$.
\begin{figure}[H]
\begin{minipage}{3.5cm}
\includegraphics[width=3.5cm]{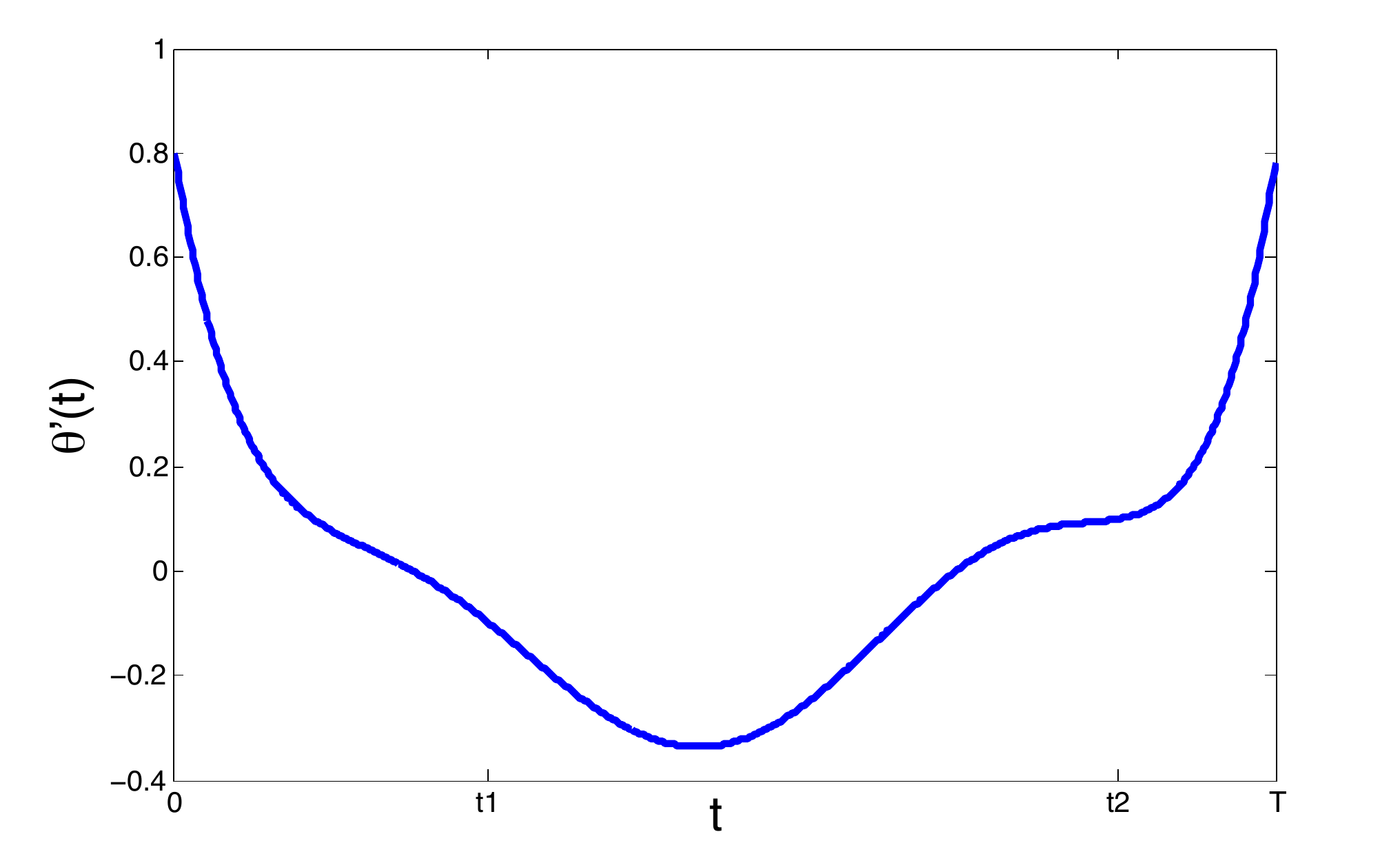}
\includegraphics[width=3.5cm]{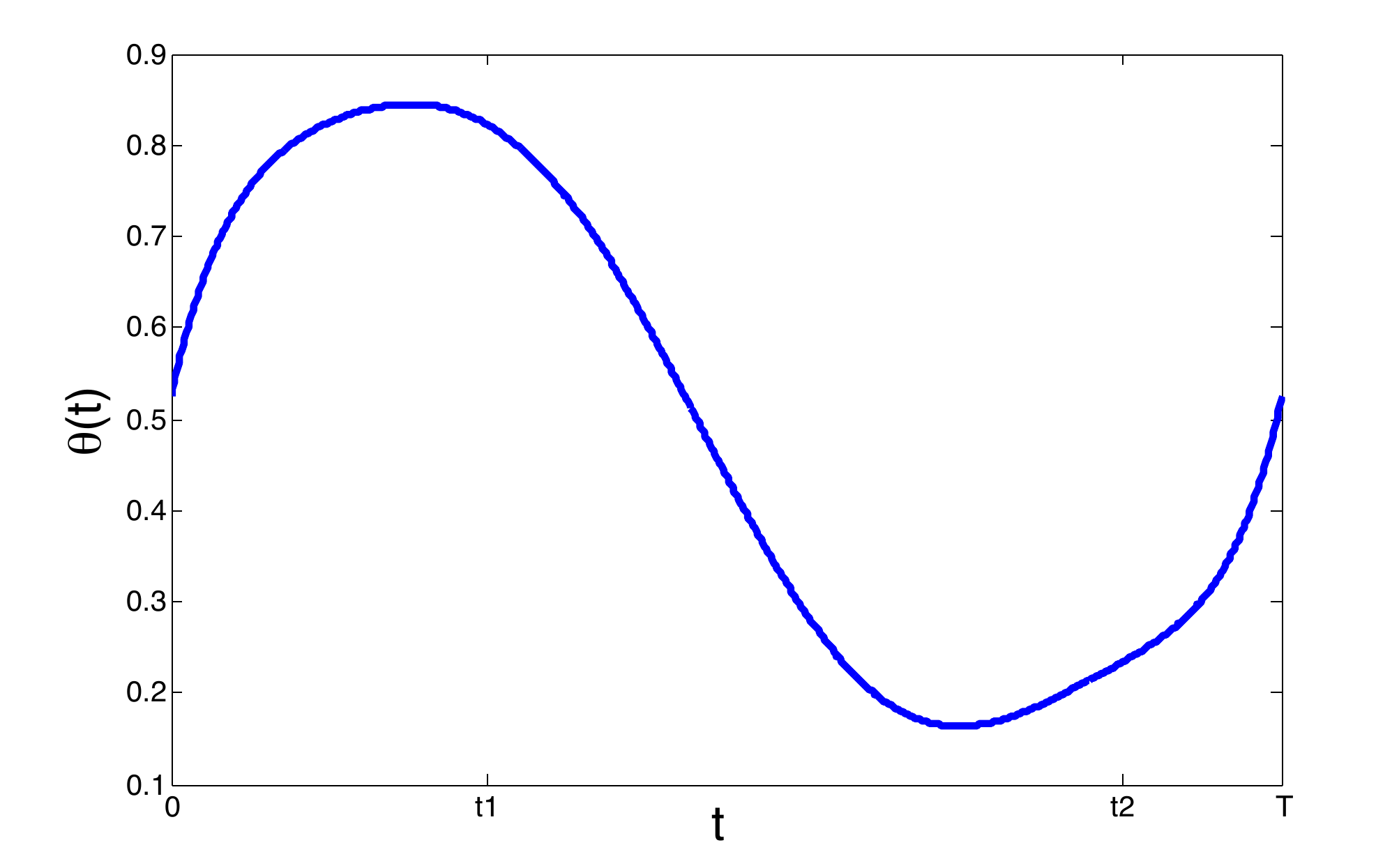}
\end{minipage}
\begin{minipage}{5cm}
\includegraphics[width=5cm]{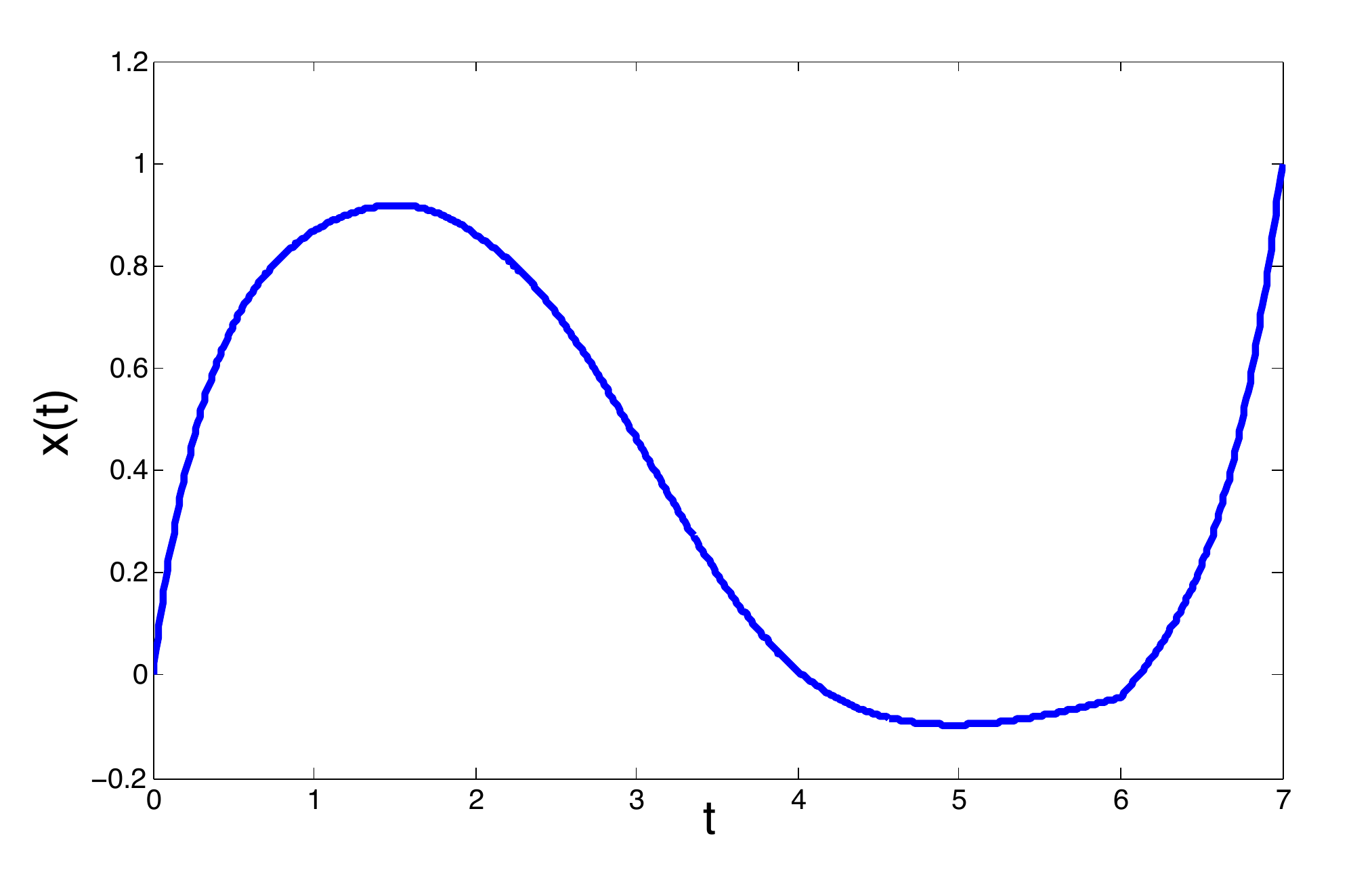}
\end{minipage}
\caption{The lpolynomial control $\dot\theta(t)$, the resulting periodic angle $\theta(t)$ and the corresponding $x$ displacement in function of time.}
\label{Theta_dot_poly}
\end{figure}

It is easy to see that (since we want a positive displacement $\theta_1>\theta_2$)   $\dot\theta$ respects the thermostatic switching rule and that after a time $T=7 \,s$ we have gained the desired displacement of $1 \,cm$.\\\\
Starting from these simulations we want to build a piecewise constant control, instead of a continuous one, to obtain the same displacement. We note that in the case of delayed thermostat a discontinuous input is in general not allowed due to the presence  of memory. The main difficulty of using a discontinuous control is to chose the switching times. Having in mind the previous simulations we can take the switching times of the continuous control and build a piecewise constant control which satisfies the constraints \eqref{constraints}.\\ Referring to the simulations in Fig. \ref{Theta_dot_poly} we get

\begin{figure}[H]
\begin{minipage}{3.5cm}
\includegraphics[width=3.5cm]{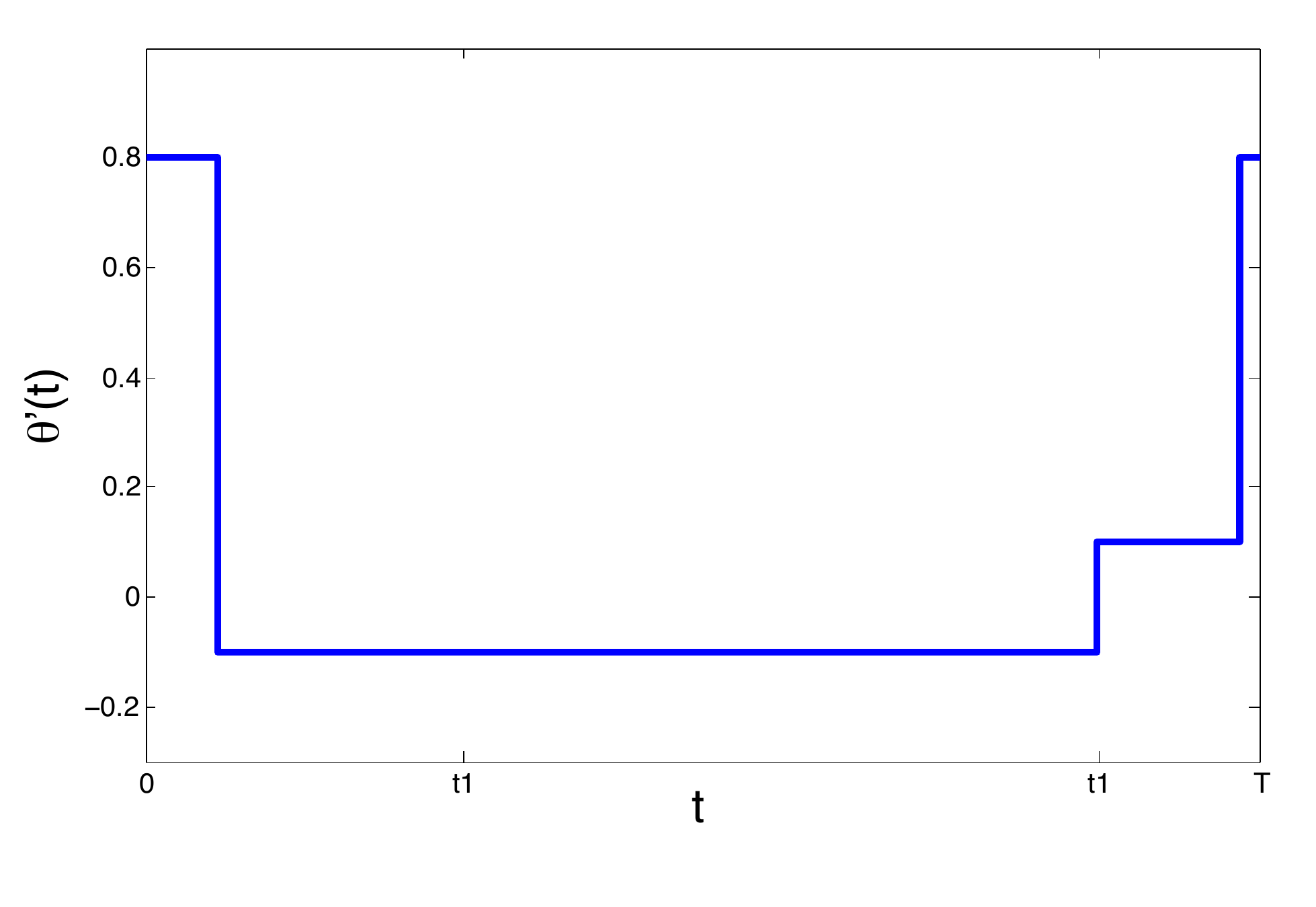}
\includegraphics[width=3.5cm]{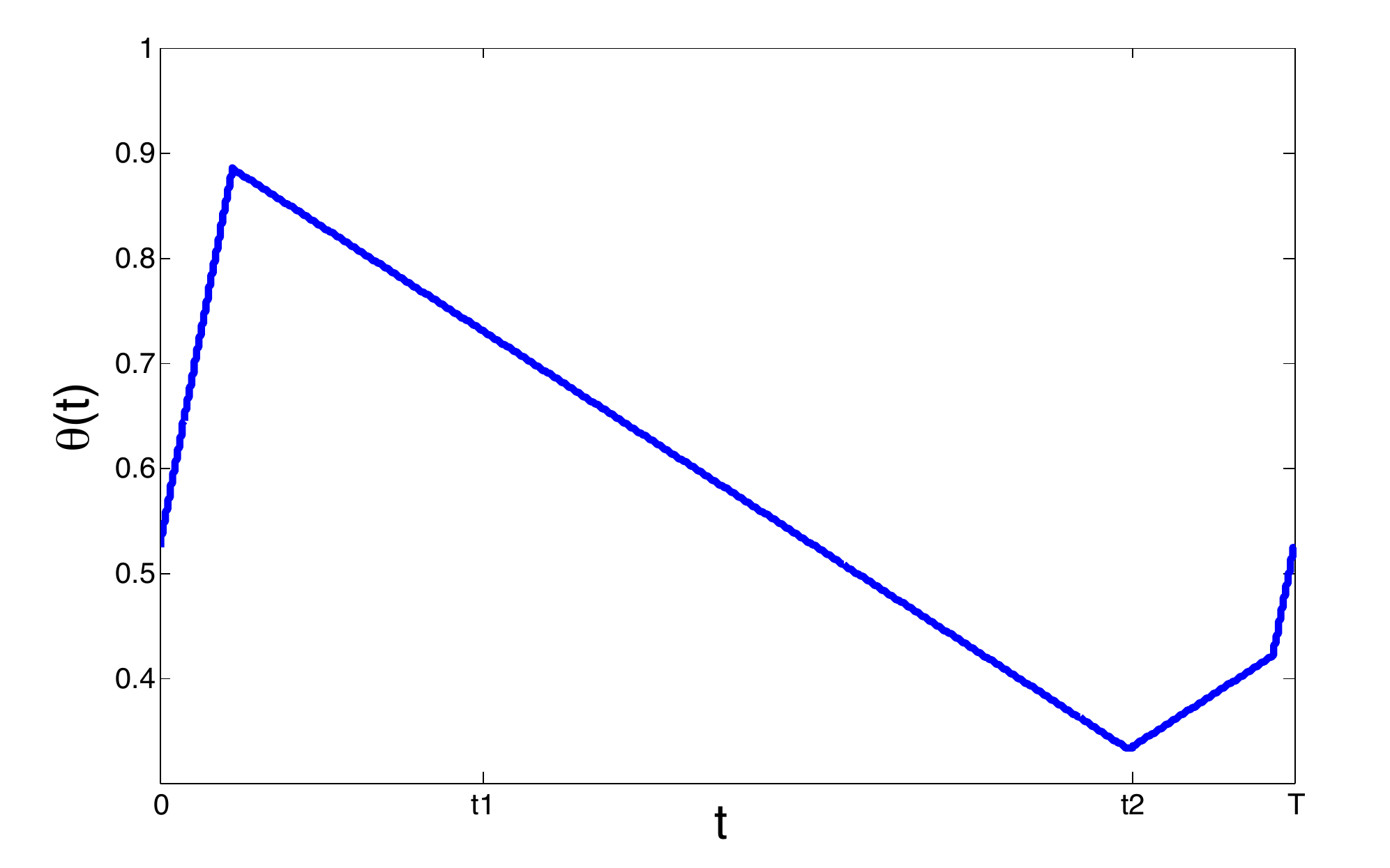}
\end{minipage}
\begin{minipage}{5cm}
\includegraphics[width=5cm]{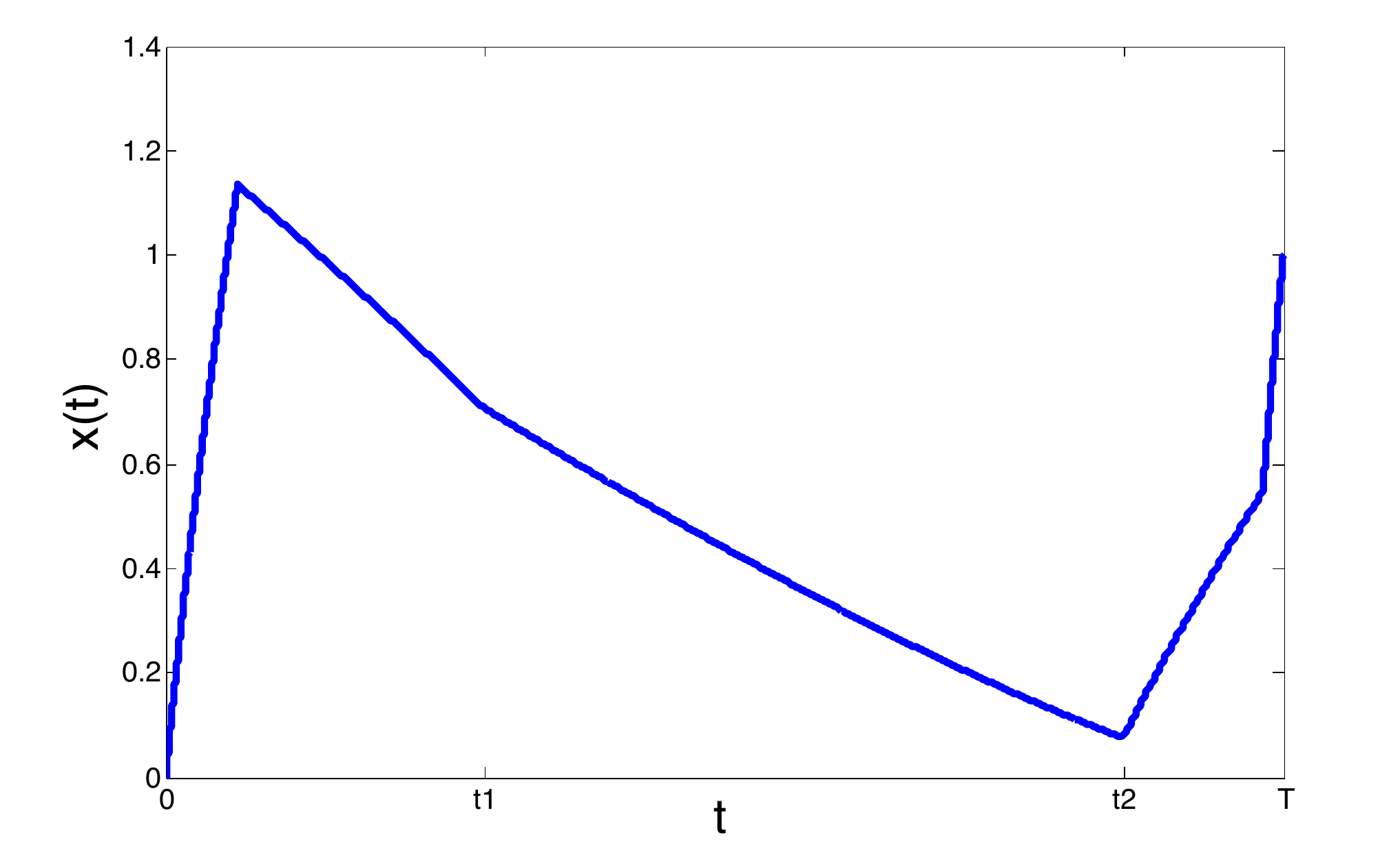}
\end{minipage}
\caption{The piecewise constant control $\dot\theta(t)$, the resulting angle $\theta$ and the corresponding $x$ displacement in function of time.}
\label{Theta_dot_linear}
\end{figure}
These simulations actually prove that the displacement does not depend on the whole control trajectory but only on the values that the angle $\theta$ and its derivative $\dot\theta$ assume in the switching times.

\section{Conclusions}
In this paper we analyze the system of a scallop proposing some strategies to overcome the famous scallop theorem. The main idea is to introduce a switching in the dynamics related to the variation of the angular velocity. This is done in two different ways, fact that helps to brake the reversibility of the equation of motion producing a net displacement. Original tool is also the introduction of the thermostat to model a delay in the change of fluid regime, and we show that it is crucial to gain both forward and backward motion. Moreover numerical simulations suggest also a way to use the switching schemes without necessarily using a continuous control input. Namely they show that it is possible to obtain the same displacement using a piecewise constant control (see Fig. \ref{Theta_dot_linear}) and choosing the switching times according to the ones used in the continuous case (see Fig. \ref{Theta_dot_poly}).

\begin{acknowledgements}
The work has been developed within the OptHySYS project of the University of Trento that is gratefully acknowledged.\\ Moreover we thank also Gruppo Nazionale Analisi Matematica Probabilit\'a e Applicazioni (GNAMPA) for partial financial support.

\end{acknowledgements}

Compliance with Ethical Standards:\\
 Funding: this study was funded by University of Trento and Gruppo Nazionale Analisi Matematica Probabilit\'a e Applicazioni (GNAMPA). \\
Conflict of Interest: the authors declare that they have no conflict of interest.


\end{document}